\documentclass[12pt]{article}
\makeatletter
\newcommand{\rmnum}[1]{\romannumeral #1}
\newcommand{\Rmnum}[1]{\expandafter\@slowromancap\romannumeral #1@}
\makeatother
\usepackage{}
\usepackage[numbers,sort&compress]{natbib}
\usepackage{color}
\usepackage{threeparttable}
\usepackage[english]{babel}
\usepackage{diagbox}
\usepackage{t1enc}
\usepackage{amsthm,amsmath,amssymb}
\usepackage{mathrsfs}
\usepackage[mathscr]{eucal}
\usepackage[latin1]{inputenc}
\usepackage[english]{babel}
\usepackage{latexsym}
\usepackage{graphicx}
\usepackage[noend]{algpseudocode}
\usepackage{algorithmicx,algorithm}
\usepackage{booktabs}
\usepackage{multirow}
\newtheorem{theorem}{Theorem}

\newtheorem{problem}{Problem}

\def \T{\textup{T}}

\begin{document}
	\title{A very short note on a problem of Godsil and Sun}
	\author{\small   Wei Wang\thanks{Email address: wangwei.math@gmail.com}
		\\
		{\footnotesize School of Mathematics, Physics and Finance, Anhui Polytechnic University, Wuhu 241000, P. R. China}
	}
	\date{}
	\maketitle
	\begin{abstract}
		Using the spectral theorem for symmetric matrices over a real closed field, we give a quick answer to a problem of Godsil and Sun on degree-similarity of graphs.		
	\end{abstract}

	

Let $G$ be a graph, and let $A=A(G)$ and $D=D(G)$ be the adjacency matrix and degree matrix of $G$, respectively. In a recent paper, Godsil and Sun \cite{godsil2025} introduced the notion of degree-similar graphs. Two graphs $G$ and $H$ are \emph{degree-similar} if there exists an invertible real matrix $M$ such that 
\begin{equation}
M^{-1}A(G)M=A(H), \text{~and~} M^{-1}D(G)M=D(H). 
\end{equation}

In \cite{wang2011}, Wang et al. defined the generalized characteristic polynomial $\psi (G,t,\mu)$ of a graph $G$ as follows:

\begin{equation}
 \psi (G,t,\mu)=\det(tI-(A-\mu D)).
\end{equation}
Clearly, if $G$ and $H$ are degree-similar, then we have  $\psi (G,t,\mu)=\psi (H,t,\mu)$.  However, the converse is not true. In \cite{durfee2015}, Durfee and Martin reported a pair of $9$-vertex graphs that have the same generalized characteristic polynomial but are not degree-similar.  Using the 16-vertex tree in an old paper of McKay \cite{mckay1977}, Godsil and Sun \cite{godsil2025} constructed infinite pairs of non-degree-similar graphs with the same generalized characteristic polynomials.  Now, we view $A-\mu D$ as a matrix over $\mathbb{Q}(\mu)$, the field of rational functions over $\mathbb{Q}$. Godsil and Sun \cite{godsil2025} proposed the following problem, which we state in a slightly different but essentially equivalent form.
\begin{problem}[\cite{godsil2025}]\label{pb}
	Let $G$ and $H$ be two graphs. Assume that $A(G)-\mu D(G)$ and $A(H)-\mu D(H)$ are similar over  $\mathbb{Q}(\mu)$. Are $G$ and $H$ degree-similar?
\end{problem}
The main result of this note is the following equivalence.
\begin{theorem}\label{main}
Let $K$ and $L$ be two symmetric matrices over $\mathbb{Q}(\mu)$. Then the following two assertions are equivalent:\\
\textup{(\rmnum{1})} $\det(tI-K)=\det(tI-L)$;\\
\textup{(\rmnum{2})} $K$ and $L$ are similar over $\mathbb{Q}(\mu)$.
\end{theorem}
Note that  matrices of the form $A-\mu D$ are symmetric matrices over $\mathbb{Q}(\mu)$. Thus, Theorem \ref{main} implies that the assumption in Problem \ref{pb} is equivalent to the assumption that $G$ and $H$ have the same characteristic polynomials. It follows that the construction proposed by Godsil and Sun \cite{godsil2025} already gives a negative answer to Problem \ref{pb}. We remark that, quite recently, using a direct but tedious argument on Smith normal forms, Fan et al.~\cite{fan2025}  proved that the construction of Godsil and Sun in some special case  satisfies the assumption of Problem \ref{pb} and hence also gives a negative answer to this problem. 

To prove Theorem \ref{main}, we recall some basic notions on formally real fields. For details, we refer to \cite{jacobson, lang}. A field $\mathbb{F}$ is said to be \emph{formally real} if $-1$ is not the sum of squares in $\mathbb{F}$. It is known that if $\mathbb{F}$ is formally real, then so is $\mathbb{F}(\mu)$ \cite[p.~634]{jacobson}. In particular, $\mathbb{Q}(\mu)$ is formally real. A field $\mathbb{F}$ is \emph{real closed} if it is formally real and no proper algebraic extension of $\mathbb{F}$ is formally real. An equivalent characterization for a field $\mathbb{F}$ to be real closed is that $\mathbb{F}$ satisfies (1) $\mathbb{F}$ is formally real; (2) the sum of two squares in $\mathbb{F}$ is a square; and (3) $\mathbb{F}(\sqrt{-1})$ is algebraically closed. Lang \cite[p.~585]{lang} remarked that the classic spectral theorem for real symmetric matrices extends naturally to any real closed field. Finally, we recall that any formally real field has an algebraic extension that is real closed. Now we can present a proof of Theorem \ref{main}.

\noindent \textbf{Proof of Theorem \ref{main}} Clearly, (\rmnum{2}) implies (\rmnum{1}). It suffices to show that (\rmnum{1}) implies (\rmnum{2}). As $\mathbb{Q}(\mu)$ is formally real, it has an extension $\mathbb{E}$ that is real closed.  Regarding $K$ and $L$ as symmetric matrices over $\mathbb{E}$, the spectral theorem implies that there exist orthogonal matrices $Q_1, Q_2$ over $\mathbb{E}$ such that  $Q_1^\T K Q_1$ and $Q_2^\T L Q_2$ are diagonal matrices $\Lambda_1$ and $\Lambda_2$. Note that the diagonal entries of $\Lambda_1$ and $\Lambda_2$ are roots of $\det(tI -K)$ and $\det(tI -L)$. Thus, the first assertion implies that $\Lambda_1$ and $\Lambda_2$ are equal (up to some permutation of diagonal entries). It follows that $K$ and $L$ are similar over $\mathbb{E}$. Since similarity of matrices over a field is invariant under field extensions, we see that  $K$ and $L$ are similar over $\mathbb{Q}(\mu)$. This completes the proof. \hfill\qedsymbol
	\section*{Acknowledgments}
	This work was partially supported by the National Natural Science Foundation of China (Grant No.~12001006) and Wuhu Science and Technology Project, China (Grant No.~2024kj015). 	
	
\end{document}